\documentclass[a4paper]{amsart}
\usepackage{amssymb} 
\usepackage[T1]{fontenc}
\usepackage[latin1]{inputenc}
\usepackage{amsfonts}
\usepackage{amsxtra}
\usepackage{ae}
\usepackage{pdfsync}
\usepackage[all]{xy}
\usepackage{enumerate}

\include{diagram}

\newcommand*{\ket}{\rangle}
\newcommand*{\bra}{\langle}

\newcommand*{\A}{\mathcal{A}}

\newcommand*{\F}{\mathcal{F}}

\renewcommand*{\H}{\mathcal{H}}

\newcommand*{\T}{\mathcal{T}}

\renewcommand*{\max}{\mathsf{f}}
\newcommand*{\red}{\mathsf{r}}

\newcommand*{\hit}{\rightharpoonup}
\newcommand*{\hitby}{\leftharpoonup}

\newcommand*{\CH}{\mathbb{C}}

\newcommand*{\LH}{\mathbb{L}}

\newcommand*{\Irr}{\mathsf{Irr}}

\DeclareMathOperator{\tr}{tr}
\DeclareMathOperator{\Tr}{Tr}

\DeclareMathOperator{\dom}{dom}
\DeclareMathOperator{\id}{id}

\newenvironment{bnum}
{\begin{list}{}
    {\setlength{\labelwidth}{15pt}
     \setlength{\leftmargin}{\labelwidth}
    }
}
{\end{list}}

\numberwithin{equation}{section}
\theoremstyle{change}
\newtheorem{theorem}{Theorem}[section]
\newtheorem{prop}[theorem]{Proposition}
\newtheorem{lemma}[theorem]{Lemma}

\newtheorem{definition}[theorem]{Definition}

\begin{document}

\title[Rapid decay]{Compact quantum metric spaces from quantum groups of rapid decay}

\author{Jyotishman Bhowmick}
\address{Indian Statistical Institute \\
203 Barrackpore Trunk Road \\
Kolkata 700 108 \\
India
}
\email{jyotishmanb@gmail.com}

\author{Christian Voigt}
\address{School of Mathematics and Statistics \\
         University of Glasgow \\
         15 University Gardens \\
         Glasgow G12 8QW \\
         United Kingdom 
}
\email{christian.voigt@glasgow.ac.uk}

\author{Joachim Zacharias}
\address{School of Mathematics and Statistics \\
         University of Glasgow \\
         15 University Gardens \\
         Glasgow G12 8QW \\
         United Kingdom 
}
\email{joachim.zacharias@glasgow.ac.uk}

\subjclass[2000]{16W30, 81R50}

\thanks{This work was supported by the Engineering and Physical Sciences Research Council Grants EP/G01419/2, EP/L013916/1 and 
the UKIERI project Quantum Probability, Noncommutative Geometry and Quantum Information. J.B. wishes to thank Sergey Neshveyev and the Department of 
Mathematics, University of Oslo, where he was a post doctoral fellow when this work started} 

\maketitle

\begin{abstract}
We present a modified version of the definition of property RD for discrete quantum groups given by Vergnioux in order to accommodate examples of
non-unimodular quantum groups. Moreover we extend the construction of spectral triples associated to discrete groups with length functions, originally 
due to Connes, to the setting of quantum groups. For quantum groups of rapid decay we study the resulting spectral triples from the point of view 
of compact quantum metric spaces in the sense of Rieffel. 
\end{abstract}

\section{Introduction} 

In the theory of noncommutative geometry in the sense of Connes \cite{Connesbook}, spectral triples can be thought of as noncommutative analogues of 
smooth Riemannian manifolds. A spectral triple $ (\A, \H, D) $ consists of a $ * $-algebra $ \A $, represented on a Hilbert space $ \H $, together with 
an unbounded self-adjoint operator $ D $ on $ \H $. The basic requirements on this data are that $ D $ has compact resolvent and that 
the commutators $ [D,a] $ are bounded for all $ a \in \A $. \\
The prototypical example of a spectral triple is given by the algebra $ \A = C^\infty(M) $ of smooth functions on a compact Riemannian spin 
manifold $ M $, acting on the Hilbert space $ \H = L^2(M, S) $ of $ L^2 $-section of the spinor bundle $ S $ of $ M $, together with the associated 
Dirac operator. Another class of examples, studied already by Connes \cite{Connescms}, arises from discrete groups equipped with length 
functions. In this case $ \A = \CH[G] $ is the complex group algebra of the group $ G $, acting on the Hilbert space $ l^2(G) $, 
and the operator $ D $ acts by multiplication with the length function. \\ 
Quantum groups can be viewed as noncommutative manifolds, and various examples of spectral triples 
have been constructed in this context, see for instance \cite{CPspectraltriple}, \cite{DLSSVdirac}, \cite{NTDirac}.
In this paper we consider more elementary examples of spectral triples for quantum groups, motivated by the construction for 
discrete groups with length functions mentioned above. Actually, the passage from discrete groups to discrete quantum groups 
is essentially straightforward in this context. Although the resulting spectral triples are trivial from the point of view of $ K $-homology, 
we show that they provide examples of quantum metric spaces in the sense of Rieffel \cite{Rieffelcqms}. In fact, in order to make a link to 
the theory of Rieffel, we have to restrict to quantum groups of rapid decay, and follow the work of Antonescu-Christensen in the group 
case \cite{ACmetricsgroup}. \\ 
The property of rapid decay for discrete quantum groups was introduced and studied by Vergnioux \cite{Vergniouxrd}, following the definition 
for classical groups in \cite{Jolissaintrd}. Quantum groups of rapid decay in the sense of \cite{Vergniouxrd} are necessarily unimodular, which unfortunately 
excludes some of the most studied examples, in particular those arising from $ q $-deformations of semisimple compact Lie groups. 
The incompatibility of the theory in \cite{Vergniouxrd} with non-unimodularity is of course invisible in the classical setting of discrete 
groups. For examples coming from $ q $-deformations it may appear somewhat surprising, because duals of classical compact Lie groups actually do have 
property RD in the sense of \cite{Vergniouxrd}. \\
In the first part of this paper we explain how a slight modification of the definitions given in \cite{Vergniouxrd} allow to remedy this situation. 
Our definitions agree with Vergnioux's for unimodular discrete quantum groups. On the other hand, we obtain a more interesting theory 
in the non-unimodular case. \\
Let us explain how the paper is organised. In section \ref{secpreliminaries} we collect some definitions from the theory 
of quantum groups and fix our notation. Section \ref{secrd} contains our modified definitions of rapid decay. In section 
\ref{secamenable} we consider amenable quantum groups and compare our notion of rapid decay with a suitable notion of polynomial growth. 
Actually, for polynomial growth the difference to the definition in \cite{Vergniouxrd} consists simply 
in replacing quantum dimensions with ordinary dimensions. In section \ref{secspectraltriples} we 
explain how the construction of a spectral triple from a group with a length function extends to the setting of quantum groups. 
The aim of section \ref{seccompactquantummetricspace} is to review the definition of compact quantum metric spaces  
in the sense of Rieffel, and to show that we obtain natural Lipschitz seminorms from the spectral triples defined in 
section \ref{secspectraltriples}. In the final section \ref{seclipnorm} we prove the Lip-norm property for suitable 
Lipschitz seminorms provided the underlying quantum group has property RD in our sense. This yields 
a family of examples of compact quantum metric spaces associated to quantum groups. \\ 
Let us make some remarks on notation. We write $ \LH(\H) $ for the space of bounded operators on a Hilbert space $ \H $. 
The closed linear span of a subset $ X $ of a Banach space is denoted by $ [X] $. Depending on the context, the symbol
$ \otimes $ denotes either the tensor product of Hilbert spaces, or the minimal tensor product of $ C^* $-algebras. 
For operators on multiple tensor products we use the leg numbering notation. 
We write $ \| \; \| = \|\; \|_{op} $ for the operator norm.

\section{Preliminaries} \label{secpreliminaries}

In this section we review some basic definitions concerning quantum groups. For more detailed information we refer to \cite{BSUM}, 
\cite{KVLCQG}, \cite{Woronowiczleshouches}. Our notation and conventions will follow \cite{NVpoincare}. \\
The main objects of study in this paper are discrete quantum groups. It is technically convenient to describe them using Hopf $ C^* $-algebras. 
Recall first that a Hopf $ C^* $-algebra is a $ C^* $-algebra $ S $ together with an injective nondegenerate $ * $-homomorphism 
$ \Delta: S \rightarrow M(S \otimes S) $, called comultiplication, such that $ (\Delta \otimes \id)\Delta = (\id \otimes \Delta) \Delta $ 
and $ [\Delta(S)(1 \otimes S)] = S \otimes S = [(S \otimes 1)\Delta(S)] $. \\ 
With this terminology, a discrete quantum group can be described by a pair of Hopf $ C^* $-algebras $ C_0(G) $ and 
$ C^*_\red(G) $ together with a multiplicative unitary $ W \in M(C_0(G) \otimes C^*_\red(G)) $, satisfying certain axioms. 
In particular, the algebra $ C^*_\red(G) $ is unital, and $ C_0(G) $ is a $ C^* $-algebraic direct sum of matrix algebras. 
We write $ \Delta $ for the comultiplication of $ C_0(G) $ and $ \hat{\Delta} $ for the comultiplication of $ C^*_\red(G) $. 
Some properties of the multiplicative unitary linking these two Hopf $ C^* $-algebras will be stated below. 
We refer to $ C_0(G) $ as the algebra of functions on $ G $, and to $ C^*_\red(G) $ as the reduced group $ C^* $-algebra of $ G $. 
The theory also provides a full group $ C^* $-algebra $ C^*_\max(G) $, which however will not show up explicitly in this paper. 
At some points we will restrict attention to the case that $ G $ is amenable, which means that the canonical quotient 
homomorphism $ C^*_\max(G) \rightarrow C^*_\red(G) $ is an isomorphism. \\ 
Inside the Hopf $ C^* $-algebra $ C_0(G) $ we have a canonical dense multiplier Hopf-$ * $-algebra $ C_c(G) $, compare \cite{vDadvances}. 
More precisely, $ C_c(G) $ is the algebraic direct sum of matrix blocks defining $ C_0(G) $. 
Moreover $ C_0(G) $ admits a left Haar weight $ \phi $, given by a positive linear functional $ \phi: C_c(G) \rightarrow \mathbb{C} $ 
satisfying $ (\id \otimes \phi)\Delta(f) = \phi(f) 1 $, and we let $ l^2(G) $ denote the GNS-construction of $ \phi $. 
We write $ \Lambda(f) \in l^2(G) $ for the image of $ f \in C_c(G) $ under the GNS-map. \\ 
The multiplicative unitary $ W $ can be considered as an element of $ \LH(l^2(G) \otimes l^2(G)) $, and we have the explicit
formula 
$$
W^*(\Lambda(f) \otimes \Lambda(g)) = (\Lambda \otimes \Lambda)(\Delta(g)(f \otimes 1)) 
$$
for its adjoint. Moreover, the $ C^* $-algebra of functions on $ G $ can be recovered from $ W $ as 
$$
C_0(G) = [(\id \otimes \LH(l^2(G))_*)(W)], 
$$
and the reduced group $ C^* $-algebra of $ G $ can be identified with  
$$
C^*_\red(G) = [(\LH(l^2(G))_* \otimes \id)(W)].   
$$
In fact, the latter formula can be taken as the definition of $ C^*_\red(G) $ if one constructs the multiplicative unitary first. \\ 
A finite dimensional unitary corepresentation of $ G $ is a unitary $ X \in C^*_\red(G) \otimes \LH(\H) $ satisfying 
$ (\hat{\Delta} \otimes \id)(X) = X_{13} X_{23} $, here $ \H $ is a finite dimensional Hilbert space, and we are using leg numbering notation. Such 
corepresentations form a semisimple $ C^* $-tensor category. We denote by $ \Irr(G) $ the set of equivalence classes of irreducible 
corepresentations of $ G $, and we write $ \epsilon $ for the trivial corepresentation on $ \mathbb{C} $. \\ 
Using corepresentation theory we can identify
$$
C_c(G) \cong \bigoplus_{\alpha \in \Irr(G)} \LH(\H_\alpha) \cong \bigoplus_{\alpha \in \Irr(G)} M_{\dim(\alpha)}(\mathbb{C})
$$
as the algebraic direct sum of the endomorphism algebras of all irreducible corepresentations. The algebra $ C_0(G) $ 
is obtained by taking the $ C^* $-algebraic direct sum instead. Finally, we will also need the algebraic multiplier algebra $ C(G) $ 
of $ C_c(G) $, which can be written as 
$$
C(G) \cong \prod_{\alpha \in \Irr(G)} \LH(\H_\alpha) \cong \prod_{\alpha \in \Irr(G)} M_{\dim(\alpha)}(\mathbb{C}), 
$$
the algebraic direct product of all endomorphism algebras $ \LH(\H_\alpha) $ for $ \alpha \in \Irr(G) $. \\ 
The matrix coefficients of all irreducible corepresentations define a canonical dense 
Hopf-$ * $-algebra $ \CH[G] \subset C^*_\red(G) $. 
The algebras $ \mathbb{C}[G] \subset C^*_\red(G) $ and $ C_c(G) \subset C_0(G) $ are linearly spanned by 
elements of the form $ (\omega \otimes \id)(W) $ and $ (\id \otimes \omega)(W) $, respectively, where 
$ \omega = \omega_{\xi, \eta} \in \LH(l^2(G))_* $ is associated to vectors $ \xi, \eta \in \Lambda(C_c(G)) $. \\ 
If $ f \in C_c(G) $ and $ x \in \mathbb{C}[G] $ are represented by $ L_f, L_x \in \LH(l^2(G))_* $ in the sense 
that $ (\id \otimes L_f)(W) = f $ and $ (L_x \otimes \id)(W) = x $, then we obtain a well-defined bilinear pairing 
\begin{equation*}
\bra f, x \ket = \bra x, f \ket = (L_x \otimes L_f)(W) = L_f(x) = L_x(f)
\end{equation*}
between $ C_c(G) $ and $ \CH[G] $, see \cite{BSUM}. \\ 
We point out that the product of $ \CH[G] $ is dual to the coproduct of $ C_c(G) $, whereas the product of $ C_c(G) $ is dual to the opposite 
coproduct of $ \CH[G] $. In other terms, we have for all $ f, g \in C_c(G)$ and $ x, y \in \CH[G] $ the relations 
\begin{equation*}
\bra f, xy \ket = \bra f_{(1)}, x \ket \bra f_{(2)}, y \ket \quad \text{and} \quad
\bra fg, x \ket = \bra f, x_{(2)} \ket \bra g, x_{(1)} \ket
\end{equation*}
where we use the Sweedler notation $ \Delta(f) = f_{(1)} \otimes f_{(2)} $ and $ \hat{\Delta}(x) = x_{(1)} \otimes x_{(2)}$
for the comultiplications on $ C_c(G)$ and $ \CH[G] $. 
Of course, this notation has to be interpreted with care, in particular, the coproduct $ \Delta(f) $ of an element $ f $ of the multiplier 
Hopf $ * $-algebra $ C_c(G) $ can be represented only as an infinite sum of simple tensors in general. \\
We shall use the notations 
\begin{align*} 
(x \hit f)(y) = f(yx), &\qquad (f \hitby x)(y) = f(xy) \\
(f \hit x)(g) = x(fg), &\qquad (x \hitby f)(g) = x(gf)
\end{align*}
for the left and right regular actions of $ \CH[G] $ on $ C_c(G) $, and of $ C(G) $ on $ \CH[G] $, respectively. Remark 
that these definitions are in accordance with our conventions for the comultiplications of $ C_c(G) $ and $ \CH[G] $. \\ 
From the duality theory of algebraic quantum groups \cite{vDadvances} it follows that there is a linear isomorphism $ \F: C_c(G) \rightarrow \CH[G] $ 
given by 
$$
\F(f)(h) = \phi(hf). 
$$
The inverse $ \F^{-1}: \CH[G] \rightarrow C_c(G) $ of this map is given by 
$$ 
\F^{-1}(x)(y) = \hat{\phi}(S(x)y) 
$$ 
where $ \hat{\phi} $ is the left and right invariant normalized Haar functional on $ \CH[G] $, and $ S $ denotes the antipode of $ \CH[G] $. 
We fix the left invariant Haar functional $ \phi $ on $ C_c(G) $ such that 
$$ 
\phi(\F^{-1}(x)) = \hat{\epsilon}(x) 
$$ 
for all $ x \in \CH[G] $, where $ \epsilon: \CH[G] \rightarrow \mathbb{C} $ denotes the counit. 
Then the map $ \F $ is isometric with respect to the standard scalar products
$$ 
\bra x, y \ket = \hat{\phi}(x^* y), \qquad \bra f, g \ket = \phi(f^*g)
$$
on $ \CH[G] $ and $ C_c(G) $, respectively. In fact, using $ \F $ we can identify $ l^2(G) $ with the GNS-representation of $ \hat{\phi} $, 
and we will write $ \hat{\Lambda}(x) \in l^2(G) $ for the image of $ x \in \CH[G] $. \\ 
The modular function for $ G $ is a multiplier $ F \in C(G) $ which relates the left and right Haar integrals of $ C_c(G) $. 
A discrete quantum group $ G $ is unimodular iff the modular function satisfies $ F = 1 $. This happens iff 
the Haar state $ \hat{\phi} $ on $ \CH[G] $ is a trace. \\
In the general case, each component $ F^\alpha $ of $ F $ for $ \alpha \in \Irr(G) $ is a positive invertible matrix such 
that $ \tr(F^\alpha) = \tr((F^\alpha)^{-1}) $, the latter being the quantum dimension of $ \alpha $, denoted by $ \dim_q(\alpha) $. One may fix bases 
such that $ F^\alpha $ is a diagonal operator for all $ \alpha \in \Irr(G) $, and we will do this in the sequel. \\
Let us write $ u^\alpha_{ij} $ for the matrix coefficients of $ \alpha \in \Irr(G) $ with respect to such an orthonormal basis 
of the representation space of $ \alpha $. Then the Schur orthogonality relations become
$$
\hat{\phi}(u^\alpha_{ij} (u^\beta_{kl})^*) = \delta_{\alpha\beta} \delta_{ik} \delta_{jl} \, \frac{(F^\alpha)_{jj}}{\dim_q(\alpha)}, 
\qquad \hat{\phi}((u^\alpha_{ij})^* u^\beta_{kl}) = \delta_{\alpha\beta} \delta_{ik} \delta_{jl} \, \frac{(F^{-1}_\alpha)_{ii}}{\dim_q(\alpha)}, 
$$
where $ \alpha, \beta \in \Irr(G) $. \\
The left and right Haar functionals for $ C_c(G) $ are given by 
$$
\phi(f) = \sum_{\alpha \in \Irr(G)} \dim_q(\alpha) \tr(F p_\alpha f), \qquad \psi(f) = \sum_{\alpha \in \Irr(G)} \dim_q(\alpha) \tr(F^{-1} p_\alpha f), 
$$
where we denote by $ p_\alpha \in C_c(G) $ the central projection corresponding to $ \alpha \in \Irr(G) $.

\section{Rapid decay for discrete quantum groups} \label{secrd}

In this section we review some definitions from \cite{Vergniouxrd} and introduce our notion of rapid decay. We also state some equivalent 
characterisations of rapid decay, following Vergnioux. \\ 
Let us first recall the notion of a length for a discrete quantum group introduced in \cite{Vergniouxrd}. 
\begin{definition} 
Let $ G $ be a discrete quantum group. A length for $ G $ is a positive element $ L \in C(G) $ such that 
\begin{bnum}
\item[a)] $ \epsilon(L) = 0 $.
\item[b)] $ S(L) = L $. 
\item[c)] $ \Delta(L) \leq L \otimes 1 + 1 \otimes L $. 
\end{bnum}
\end{definition}
These conditions reflect the classical definition of length functions. Note that $ L $ can be viewed as a sequence of positive matrices indexed 
by the elements of $ \Irr(G) $. \\
We will be mainly interested in central lengths, that is, lengths $ L $ which are central elements 
of the algebra $ C(G) $. Such central lengths are obtained from length functions $ l $ on $ \Irr(G) $ in the following sense. 
\begin{definition} 
Let $ G $ be a discrete quantum group. A length function for $ G $ is a function $ l: \Irr(G) \rightarrow [0, \infty) $ such that 
\begin{bnum}
\item[a)] $ l(\epsilon) = 0 $ where $ \epsilon \in \Irr(G) $ is the trivial corepresentation. 
\item[b)] $ l(\overline{\alpha}) = l(\alpha) $ for all $ \alpha \in \Irr(G) $. 
\item[c)] $ \alpha \subset \beta \otimes \gamma $ implies $ l(\alpha) \leq l(\beta) + l(\gamma) $ for all $ \alpha, \beta, \gamma \in \Irr(G) $. 
\end{bnum}
The length function $ l $ is called proper if for any $ n \in \mathbb{N} $ there are 
only finitely many irreducible corepresentations $ \alpha \in \Irr(G) $ with $ l(\alpha) \leq n $ and $ l(\alpha) = 0 $ iff $ \alpha = \epsilon $. 
\end{definition} 
A length function $ l $ induces a central length $ L \in C(G) $ by the formula 
$$
L = \sum_{\alpha \in \Irr(G)} l(\alpha) p_\alpha, 
$$ 
recall that $ p_\alpha \in C_c(G) $ is the unit element in the matrix block corresponding to $ \alpha $. We can also view $ L $ as an unbounded 
self-adjoint operator on $ l^2(G) $ in the obvious way. Any central length arises in this way, and we will freely pass from $ l $ 
to $ L $ in the sequel. \\ 
Basic examples of length functions are given by word length functions on finitely generated quantum groups. We recall that a subset $ D \subset \Irr(G) $ 
is said to generate the discrete quantum group $ G $ iff every corepresentation $ \alpha \in \Irr(G) $ is contained in some iterated tensor product 
of corepresentations from $ D $. The quantum group $ G $ is called finitely generated provided there exists a finite subset $ D \subset \Irr(G) $ which 
generates $ G $. \\ 
Given a finitely generated quantum group $ G $ with finite generating set $ D $, we obtain a proper length function $ l_D $ on $ \Irr(G) $ by 
letting $ l_D(\alpha) $ be the smallest number $ k $ such that $ \alpha \subset \alpha_1 \otimes \cdots \otimes \alpha_k $ and $ \alpha_j \in D $ for all $ j $. 
Although $ l_D $ clearly depends on $ D $, it can be shown that the definitions and results in the sequel 
do not depend on the choice of the generating set in an essential way, compare lemma 3.3 and remark 3.6 in \cite{Vergniouxrd}. \\
Given a length function $ l $ on $ G $, we let $ p_n \in \LH(l^2(G)) $ be the sum of all projections $ p_\alpha $ for 
$ \alpha \in \Irr(G) $ such that $ \|L p_\alpha \|_{op} \in (n - 1, n] $. Notice that $ p_n $ is a finite rank projection for all $ n \in \mathbb{N}_0 $ 
iff $ l $ is proper. \\
Let $ G $ be a discrete quantum group. In the sequel we shall work with the self-adjoint element $ C $ of $ C(G) $ given by 
$$
C = \sum_{\alpha \in \Irr(G)} \frac{\dim_q(\alpha)}{\dim(\alpha)} F^\alpha p_\alpha. 
$$ 
We remark that $ C = 1 $ iff $ G $ is unimodular. 
If $ L $ is a length on $ G $, we define the associated Sobolev $ s $-norm for $ s \geq 0 $ and $ f \in C_c(G) $ by 
\begin{align*}
\| f \|_{2,s}^2 &= \sum_{\alpha \in \Irr(G)} \frac{\dim_q(\alpha)^2}{\dim(\alpha)} \tr(p_\alpha ((1 + L)^s f F)^* (1 + L)^s f F) \\
&= \phi(((1 + L)^s f C^{1/2})^* (1 + L)^s f C^{1/2}) \\
&= \bra (1 + L)^s f C^{1/2}, (1 + L)^s f C^{1/2} \ket. 
\end{align*}
Observe that $ \|\; \|_{2,0} = \|\;\|_2 $ iff $ G $ is unimodular. 
Similarly, for $ s \geq 0 $ and $ x \in \CH[G] $ we define
\begin{align*}
\| x \|_{2,s}^2 &= \bra (1 + L)^s \F^{-1}(C^{1/2} \hit x), (1 + L)^s \F^{-1}(C^{1/2} \hit x) \ket \\
&= \sum_{\alpha \in \Irr(G)} \frac{\dim_q(\alpha)^2}{\dim(\alpha)} \tr(p_\alpha ((1 + L)^s \F^{-1}(F \hit x))^* (1 + L)^s \F^{-1}(F \hit x)), 
\end{align*}
where we recall that $ \F^{-1}: \CH[G] \rightarrow C_c(G) $ is the Fourier transform given by $ \F^{-1}(x)(y) = \hat{\phi}(S(x)y) $ 
for $ x, y \in \mathbb{C}[G] $. Notice that 
$$ 
\F(f C^{1/2})(h) = \phi(h f C^{1/2}) = \phi(C^{1/2} h f) = (C^{1/2} \hit \F(f))(h), 
$$
which implies 
$$
\F^{-1}(C^{1/2} \hit x) = \F^{-1}(x) C^{1/2}. 
$$
Hence $ \F $ is an isometry with respect to the $ \| \; \|_{2,s} $-norms. 
\begin{definition} \label{defrd}
Let $ G $ be a discrete quantum group. We say that $ G $ has property RD with respect to a central length $ L $ on $ G $ if there exist 
constants $ c, s > 0 $ such that 
$$ 
\| \F(f) \|_{op} \leq c \| f \|_{2,s} 
$$ 
for all $ f \in C_c(G) $. \\
We say that $ G $ has property RD if it has property RD with respect to some central length. 
\end{definition} 
Since $ C = 1 $ for a unimodular discrete quantum group, our definition 
of the Sobolev norms above reduces to the definitions in \cite{Vergniouxrd} in this case. In other words, for 
unimodular quantum groups definition \ref{defrd} is equivalent to the definition of property RD given by Vergnioux. \\ 
As we will see below, this is not true for non-unimodular quantum groups. One might therefore call the property defined above 
modular property RD, or something alike, in order to make a distinction with the original notion. However, since there 
is no real conflict in terminology neither in the unimodular nor in the non-unimodular case - a non-unimodular quantum group 
can only possibly have property RD according to the above version of the definition - we have refrained from introducing new
terminology. \\ 
Our above definition of Sobolev norms may appear somewhat arbitrary, but it can be motivated as follows. Essentially, the 
operator $ C^{1/2} $ controls the deviation of the norm 
$$ 
\|\hat{\Lambda}(u^\alpha_{ij})\| = \frac{(F^\alpha)^{-1/2}_{ii}}{\dim_q(\alpha)^{1/2}} 
$$
from its value $ \dim(\alpha)^{-1/2} $ in the unimodular case. More precisely, we have
$$
\|\hat{\Lambda}(C^{1/2} \hit u^\alpha_{ij})\| = \frac{1}{\dim(\alpha)^{1/2}}, 
$$
so that the action of $ C^{1/2} $ compensates for the rescaling of Hilbert space norms. We remark that one could also 
incorporate $ C^{1/2} $ in the definition of the Fourier transform, which would simplify some formulas in the sequel. 
The remaining ingredients in the definition of the Sobolev norms are as in \cite{Vergniouxrd}, with the only difference that 
we are working with left Haar weights. \\ 
If $ L $ is a length on $ G $ we denote by $ \H^s_L(G) $ the completion of $ \mathbb{C}[G] $ with respect to the 
Sobolev $ s $-norm $ \|\; \|_{2,s} $. Moreover we define the associated Schwartz space $ \H^\infty_L(G) $ by 
$$
\H^\infty_L(G) = \bigcap_{s \geq 0} \H^s_L(G). 
$$
The space $ \H^\infty_L(G) $ is naturally a Fr\'echet space with respect to the topology given by the Sobolev seminorms. 
Let us prove the following variant of proposition 3.5 in \cite{Vergniouxrd}. 
\begin{prop} \label{RDchar}
Let $ G $ be a discrete quantum group and let $ L $ be a central length on $ G $. Then the following conditions are 
equivalent. 
\begin{bnum}
\item[a)] $ G $ has property RD with respect to $ L $. 
\item[b)] There exists $ c, s > 0 $ such that $ \| x \|_{op} \leq c \| x \|_{2,s} $ for all $ x \in \CH[G] $. 
\item[c)] The identity map $ \CH[G] \rightarrow \CH[G] $ induces a continuous linear embedding
$$
\H^\infty_L(G) \subset C^*_\red(G).
$$ 
\item[d)] There exists a polynomial $ p(x) \in \mathbb{R}[x] $ such that 
$$ 
\| \F(f)\|_{op} \leq p(n) \| f \|_{2,0} 
$$ 
for all $ n \in \mathbb{N} $ and $ f \in p_n C_c(G) $. 
\item[e)] There exists a polynomial $ p(x) \in \mathbb{R}[x] $ such that 
$$ 
\|p_l \F(f)p_k\|_{op} \leq p(n) \| f \|_{2,0} 
$$ 
for all $ k,l,n \in \mathbb{N} $ and $ f \in p_n C_c(G) $. 
\end{bnum}
\end{prop}
\proof The arguments follow precisely the pattern of \cite{Vergniouxrd}, for the convenience of the reader we shall include the details. \\ 
$ a) \Leftrightarrow b) $ This follows immediately from the fact that the Fourier transform $ \F: C_c(G) \rightarrow \CH[G] $ 
is a linear isomorphism preserving Sobolev norms. \\
$ b) \Rightarrow c) $ From the norm estimate in $ b) $ we obtain the existence of a continuous linear map $ \iota: \H^\infty_L(G) \rightarrow C^*_\red(G) $, 
and we may compose with the embedding $ C^*_\red(G) \rightarrow l^2(G) \subset C(G) $. Since the completion of $ C_c(G) $ 
with respect to $ \|\; \|_s $ can be realised as a subspace of $ C(G) $ in a compatible way it follows that $ \iota $ is injective. \\
$ c) \Rightarrow b) $ is obvious. \\ 
$ a) \Rightarrow d) $ Let $ s $ be the natural number  and $ c > 0 $ such that $ \| \F(f) \|_{op} \leq c \| f \|_{2,s} $. Then for $ f \in p_n C_c(G) $ 
we have
\begin{align*}
\|\F(f)\|_{op} \leq c\|f\|_{2,s} = c\|(1 + L)^s f \|_{2,0} \leq c\|(1 + n)^s f \|_{2,0} \leq p(n) \|f\|_{2,0}
\end{align*}
where $ p(x) = c(1 + x)^s $. \\
$ d) \Rightarrow a) $ Let us choose constants $ c_1, s > 0 $ such that $ p(n) \leq c_1 (1 + n)^s $ for all $ n \in \mathbb{N} $. Then we obtain 
\begin{align*}
\|\F(f) \|_{op} &\leq \sum_{n = 0}^\infty \|\F(p_n f)\|_{op} \\
&\leq c_1 \sum_{n = 0}^\infty (1 + n)^s \| p_n f \|_{2,0} \\
&= c_1 \sum_{n = 0}^\infty \frac{1}{1 + n} (1 + n)^{s + 1} \| p_n f \|_{2,0} \\
&\leq c_1 \biggl(\sum_{n = 0}^\infty \frac{1}{(1 + n)^2}\biggr)^{1/2} \biggl(\sum_{n = 0}^\infty (1 + n)^{2s + 2} \| p_n f \|_{2,0}^2 \biggr)^{1/2} \\
&\leq c_2 \|(1 + L)^{s + 1} f\|_{2,0} = c_2 \|f\|_{2,s + 1} 
\end{align*}
for a suitable constant $ c_2 $, using the Cauchy-Schwarz inequality as in the classical case \cite{Haagerupnonnuclear}. \\ 
$ d) \Rightarrow e) $ is obvious. \\ 
$ e) \Rightarrow d) $ Let $ f \in p_n C_c(G) $. Then if $ p_l \F(f) p_k \neq 0 $ we have $ (k,l,n) \in \T $ where 
$ \T \subset \mathbb{N}^3 $ is the set of all $ (k,l,n) $ such that $ \Delta(p_n)(p_k \otimes p_l) \neq 0 $. Indeed, 
due to lemma 3.4 in \cite{Vergniouxrd} the set $ \T $ is stable under permutations and we have 
$$
p_l \F(f) p_k = (\phi \otimes \id)((1 \otimes p_l) W (f \otimes p_k)) = (\phi \otimes \id)(W \Delta(p_l)(p_n f \otimes p_k)), 
$$
using $ \F(f) \Lambda(h) = \bra S^{-1}(h_{(1)}), \F(f) \ket \Lambda(h_{(2)}) = \phi(S^{-1}(h_{(1)}) f) \Lambda(h_{(2)}) $ for 
$ h \in C_c(G) $, and $ (\phi \otimes \id)(W (f \otimes 1)) \Lambda(h) = \phi(S^{-1}(h_{(1)})f) \Lambda(h_{(2)}) $. \\ 
For a unit vector $ \xi \in l^2(G) $ we obtain 
\begin{align*}
\|\F(f) \xi \|^2 &= \sum_{l = 0}^\infty \|p_l \F(f) \xi\|^2 \\
&\leq \sum_{l = 0}^\infty \biggl(\sum_{k = 0}^\infty \|p_l \F(f) p_k \xi\| \biggr)^2 \\
&\leq p(n)^2 \|f\|_{2,0}^2 \sum_{l = 0}^\infty \biggl(\sum_{k \mid (k,l,n) \in \T} \|p_k \xi \| \biggr)^2
\end{align*}
using condition $ e) $. 
According to lemma 3.4 in \cite{Vergniouxrd}, for fixed $ l \in \mathbb{N} $ the number of elements $ k $ such 
that $ (k,l,n) \in \T $ is bounded above by $ 2n + 5 $. Hence by the Cauchy-Schwarz inequality we have 
$$
\biggl(\sum_{k \mid (k,l,n) \in \T} \|p_k \xi \|\biggr)^2 \leq (2n + 5) \biggl(\sum_{k \mid (k,l,n) \in \T} \|p_k \xi \|^2 \biggr)
$$ 
for any $ l $. Therefore, 
\begin{align*}
\sum_{l = 0}^\infty \biggl(\sum_{k \mid (k,l,n) \in \T} \|p_k \xi \|\biggr)^2 
&\leq \sum_{l = 0}^\infty (2n + 5) \biggl(\sum_{k \mid (k,l,n) \in \T} \|p_k \xi \|^2 \biggr) \\
&\leq (2n + 5)^2 \sum_{k = 0}^\infty \|p_k \xi \|^2 = (2n + 5)^2 \|\xi\|^2, 
\end{align*}
and we get $ \|\F(f)\|_{op} \leq (2n + 5) p(n) \|f\|_{2,0} $ as desired. \qed

\section{Property RD and polynomial growth} \label{secamenable} 

In this section we study property RD for amenable quantum groups, and we show that it is equivalent to a suitable notion of 
polynomial growth, again following the work of Vergnioux. However, we have to modify the definition of polynomial growth 
introduced in \cite{Vergniouxrd}. We will comment on the relation between the various concepts below. \\
Let us start with the following definition. 
\begin{definition} \label{defpolynomialgrowth}
Let $ G $ be a discrete quantum group. We say that $ G $ has polynomial growth with respect to a central length function $ l $ on $ G $ if there 
exists a polynomial $ p(x) \in \mathbb{R}[x] $ such that 
$$
\sum_{l(\alpha) \in (n - 1, n]} \dim(\alpha)^2 \leq p(n) 
$$
for all $ n \in \mathbb{N} $. \\
We say that $ G $ has polynomial growth if it has polynomial growth with respect to some central length function on $ G $. 
\end{definition} 
A finitely generated discrete quantum group $ G $ has polynomial growth iff it has polynomial growth with respect to a word length 
function $ l $ on $ G $. \\
If $ l $ is a central length function on $ G $ we shall write $ S^n \subset \Irr(G) $ for the set of all corepresentations $ \alpha $ 
satisfying $ l(\alpha) \in (n - 1, n] $. We note that $ G $ has polynomial growth with respect to $ l $ iff the sequences 
$$
s_n = |S^n|, \qquad d_n = \sup_{\alpha \in S^n} \dim(\alpha)^2
$$ 
both have polynomial growth. In particular, a unimodular discrete quantum group has polynomial growth in the sense of \cite{Vergniouxrd} iff 
it has polynomial growth in the sense of definition \ref{defpolynomialgrowth}. \\ 
For non-unimodular quantum groups the concept introduced above differs from the notion in \cite{Vergniouxrd}. Indeed, for polynomial growth 
in the sense of Vergnioux one has to replace classical dimensions by quantum dimensions in definition \ref{defpolynomialgrowth}. 
In order to distinguish the two notions one could refer to them as classical and quantum polynomial growth, respectively. As in the 
case of property RD we shall however refrain from using new terminology. 
We observe that a non-unimodular quantum group cannot have polynomial growth in the sense of \cite{Vergniouxrd}, so this should not lead 
to confusion. \\ 
We have the following version of a result in \cite{Vergniouxrd}.  
\begin{prop} \label{rdpgequiv} 
Let $ G $ be a discrete quantum group and let $ l $ be a length function on $ G $. 
\begin{bnum} 
\item[a)] If $ G $ has polynomial growth with respect to $ l $ then $ G $ has property RD with respect to $ l $. 
\item[b)] If $ G $ is amenable and has property RD with respect to $ l $, then $ G $ has polynomial growth with respect to $ l $. 
\end{bnum}
In particular, rapid decay and polynomial growth are equivalent in the amenable case. 
\end{prop} 
\proof Again, we follow the arguments in \cite{Vergniouxrd}, and include the details for the convenience of the reader. \\ 
$ a) $ Recall that every element of $ \CH[G] \subset C^*_\red(G) $ can be written in the 
form $ (\omega \otimes \id)(W) $ where $ W \in M(C_0(G) \otimes C^*_\red(G)) $ is the multiplicative unitary of $ G $ and 
$ \omega $ a linear functional on $ C_0(G) $ of the form $ \omega = \F(f) $ for $ f \in C_c(G) $. Here $ \F $ denotes the Fourier transform as 
above. Notice that $ (\omega \otimes \id)(W) = \F(f) \in C^*_\red(G) $ for $ \omega = \F(f) $, and 
$ \|\F(f)\|_{op} = \|(\omega \otimes \id)(W) \| \leq \| \omega \| $. \\  
Let us assume first that $ f \in p_n C_c(G) $ is such that $ \omega = \F(f) $ is a positive linear functional. The latter is equivalent to saying 
that $ f F \in C_c(G) $ is positive. Taking an approximate identity $ (u_j)_{j \in J} $ for $ C_0(G) $ of central projections we 
compute $ \|\omega\| = \lim \F(f)(u_j) = \lim \phi(u_j f) = \phi(f) $. Hence we obtain
\begin{align*}
\| \F(f) \|_{op} &\leq \|\omega\| \\ 
&= \phi(f) \\ 
&= \phi(p_n f) \\
&= \sum_{\alpha \in S^n} \dim_q(\alpha) \tr(p_\alpha f F) \\
&= \sum_{\alpha \in S^n} \dim_q(\alpha)^{1/2} \dim(\alpha)^{1/2} \tr(p_\alpha f F^{1/2} C^{1/2}) \\
&\leq \biggl(\sum_{\alpha \in S^n} \tr(\dim(\alpha) p_\alpha)\biggr)^{\frac{1}{2}}
\biggl(\sum_{\alpha \in S^n} \dim_q(\alpha) \tr(p_\alpha F^{\frac{1}{2}} (f C^{\frac{1}{2}})^* f C^{\frac{1}{2}} F^{\frac{1}{2}})\biggr)^{\frac{1}{2}} \\
&\leq \biggl(\sum_{\alpha \in S^n} \dim(\alpha)^2 \biggr)^{1/2}\|f\|_{2,0} \\
&\leq \sqrt{p(n)} \, \|f\|_{2,0} 
\end{align*}
by the Cauchy-Schwarz inequality, where $ p(x) $ is the polynomial appearing in the polynomial growth estimate for $ G $. \\ 
Next assume that the functional $ \omega = \F(f) $ is hermitian, which is equivalent to saying that $ f F $ is self-adjoint. In this case 
we may write $ f = f_+ - f_- $ where $ f_\pm F^{-1} $ are positive. In fact, we have $ f_\pm = \pm e_\pm f $ for suitable projections $ e_\pm $. The previous computation yields
\begin{align*}
\| \F(f) \|_{op} &\leq \| \F(f_+) \|_{op} + \| \F(f_-) \|_{op} \\
&\leq \sqrt{p(n)}(\|e_+ f\|_{2,0} + \|e_- f\|_{2,0}) \leq \sqrt{2 p(n)} \|f\|_{2,0}, 
\end{align*}
where the last step follows by inspecting the definition of the norm $ \|\; \|_{2,0} $. \\
Finally, consider an arbitrary element $ f \in p_n C_c(G) $ and write $ f $ as sum $ f = g + ih $ such that both $ gF $ and $ hF $ 
are self-adjoint. In this case, taking into account our previous computations, we obtain 
\begin{align*}
\| \F(f) \|_{op} &\leq \| \F(g) \|_{op} + \| \F(h) \|_{op} \\
&\leq \sqrt{2p(n)} (\|g\|_{2,0} + \|h\|_{2,0}) \\
&\leq \sqrt{4 p(n)} \|f\|_{2,0} 
\end{align*}
using that 
$$ 
\|g\|_{2,0}^2 + \|h\|^2_{2,0} = \|g + ih\|_{2,0}^2 = \|f\|_{2,0}^2 
$$ 
since 
$$ 
\tr(p_\alpha (gF)^* hF) = \tr(p_\alpha gF hF) = \tr(p_\alpha(hF)^* gF) 
$$ 
for any $ \alpha \in \Irr(G) $. It follows that $ G $ has property RD. \\
$ b) $ Due to amenability, the counit $ \hat{\epsilon}: \CH[G] \rightarrow \CH $ extends to a $ * $-homomorphism $ C^*_\red(G) \rightarrow \CH $. Moreover, 
we have $ \hat{\epsilon}(\F(f)) = \phi(f) $ for all $ f \in C_c(G) $. Therefore 
$$
|\phi(f)| = |\hat{\epsilon}(\F(f))| \leq \| \F(f)\|_{op} \leq c \|f \|_{2,s} 
$$
for some constants $ c, s > 0 $. Let us define
$$ 
q_n = \sum_{\alpha \in \Irr(G)} \frac{\dim(\alpha)}{\dim_q(\alpha)} p_\alpha p_n F^{-1} 
= \sum_{\alpha \in S^n} \frac{\dim(\alpha)}{\dim_q(\alpha)} p_\alpha F^{-1} 
= \sum_{\alpha \in S^n} p_\alpha C^{-1}, 
$$ 
and observe that 
$$ 
\phi(q_n) = \sum_{\alpha \in S^n} \dim(\alpha)^2. 
$$
Hence $ q_n $ can be used to test polynomial growth. More precisely, we compute
\begin{align*}
\phi(q_n)^2 &\leq c^2 \|(1 + L)^s q_n\|_{2,0}^2 \\ 
&\leq c^2 (1 + n)^{2s} \sum_{\alpha \in S^n} \frac{\dim_q(\alpha)^2}{\dim(\alpha)} \tr(p_\alpha C^{-2} F^2) \\
&= c^2 (1 + n)^{2s} \sum_{\alpha \in S^n} \dim_q(\alpha) \tr(p_\alpha C^{-1} F) \\
&= c^2 (1 + n)^{2s} \phi(q_n) 
\end{align*}
which implies $ \phi(q_n) \leq c^2 (1 + n)^{2s} $. That is, $ G $ has polynomial growth. \qed \\ 
Using proposition \ref{rdpgequiv} we conclude that all duals of $ q $-deformations of compact semisimple Lie groups have property RD. 
Indeed, these discrete quantum groups are amenable \cite{Banicacompactsubfactor}, and the Weyl dimension formula 
implies that they have polynomial growth in the sense of definition \ref{defpolynomialgrowth}, 
see example 4.5 in \cite{Vergniouxrd}. We record the following precise statement of this fact. 
\begin{prop} 
Let $ q \in (0,1] $ and let $ G_q $ be the standard deformation of a simply connected compact semisimple Lie group $ G $. Then the discrete dual 
quantum group $ \hat{G}_q $ has property RD. 
\end{prop} 
This shows in particular that, in the non-unimodular case, our definition of property RD differs from the original definition given by Vergnioux.

\section{Spectral triples from length functions} \label{secspectraltriples}

In this section we explain how to associate spectral triples to discrete quantum groups equipped with proper length functions. Moreover we 
study basic properties of the spectral triples obtained this way. \\
Let us first recall the definition of a spectral triple due to Connes, see \cite{Connesbook}, \cite{GFVbook}. 
\begin{definition} 
A spectral triple $ (\A, \H, D) $ consists of 
\begin{bnum} 
\item[a)] a $ * $-algebra $ \A $, faithfully represented on $ \H. $
\item[b)] a (graded) Hilbert space $ \H $, together with 
\item[c)] an unbounded self-adjoint (odd) operator $ D $ in $ \H $ 
\end{bnum} 
such that 
\begin{bnum} 
\item[a)] the commutators $ [D, a] $ are bounded for all $ a \in \A $
\item[b)] $ D $ has compact resolvent, that is, 
$$
(1 + D^2)^{-1}
$$
is compact. 
\end{bnum}
\end{definition} 
In the above definition we have indicated how to include a grading in the general setup, but this will not be relevant in our examples. \\
The prototypical example of a spectral triple is given by $ \A = C^\infty(M) $, the algebra of smooth functions on a compact 
Riemannian spin manifold, the Hilbert space $ \H = L^2(M, S) $ of $ L^2 $-sections of the spinor bundle of $ M $, and the 
Dirac operator $ D $. In this case the formula 
$$ 
d(x,y) = \sup_{f \in C^\infty(M) \mid \| [D, f] \| \leq 1}  |f(x) - f(y)| 
$$
allows to express the Riemannian metric, or rather the corresponding geodesic distance function, in terms of operator theoretic data. \\
Connes studied spectral triples associated to discrete groups with length functions \cite{Connescms}. Let $ G $ be a discrete group and let $ l $ be a 
proper length function on $ G $. Consider the Hilbert space $ l^2(G) $ and the unbounded operator $ D $ in $ l^2(G) $ defined 
on $ \CH[G] \subset l^2(G) $ by 
$$
D = \sum_{n \in \mathbb{N}_0} n p_n
$$
where $ p_n = q_n - q_{n - 1} $, and $ q_k $ denotes the projection onto the finite dimensional subspace of $ l^2(G) $ spanned by group elements of 
length at most $ k $. If $ l $ comes from a word metric we may view $ D $ as multiplication by the length function $ l $. \\ 
It is not hard to check that one obtains a spectral triple in this way. In fact, we will consider a more general class of spectral 
triples below, and for this it is convenient to work within the framework of filtered algebras. More precisely, 
let $ \A $ be a complex unital dense $ * $-subalgebra of a unital $ C^* $-algebra $ A $. 
Let us call $ \A $ a filtered algebra if we are given a family of finite dimensional subspaces $ (\A_m)_{m \in \mathbb{N}_0} $ of $ \A $ such that
\begin{bnum}
\item[a)] $ \A = \bigcup^{\infty}_{n = 0} \A_n $, 
\item[b)] $ \A_0 = \mathbb{C} 1 $, 
\item[c)] $ \A_m \subseteq \A_n $ for $ m < n $, 
\item[d)] $ \A^*_n = \A_n $ for all $ n $, 
\item[e)] $ \A_m \A_n \subseteq \A_{m + n} $ for all $ m,n $. 
\end{bnum}
Assume further that $ \sigma $ is a state on $ A $ which is faithful on $ \A $, and let $ \H = L^2(A, \sigma) $ be the GNS-construction 
of $ \sigma $. Then each $ \A_n $ is a finite dimensional, hence closed, subspace of $ \H $. 
Let $ q_n $ denote the orthogonal projection of $ \H $ onto $ \A_n $. 
Let $ p_n = q_n - q_{n - 1}  $ and $ p_0 = q_0 $. Then the formula 
$$ 
D = \sum_{n \in \mathbb{N}_0} n p_n
$$ 
defines an unbounded operator $ D $ on $ \H $. \\ 
Filtrations and associated Dirac operators were first introduced by Voiculescu in \cite{Voiculescuquasicentral} and further studied in \cite{OzawaRieffel}. 
The following general fact is a restatement of lemma 1.1. in \cite{OzawaRieffel}. 
\begin{lemma} \label{ozawarieffelsptriple}
The triple $ (\A, \H, D) $ associated to a filtered algebra $ \A $ together with a state as above is a spectral triple. 
\end{lemma}
We note that the $ K $-homology class of the resulting (odd) triple is trivial because the operator $ D $ is positive. That is, the 
corresponding phase $ F $ is identically $ 1 $. In general, one should rather think of $ D $ as 
the absolute value $ |D| $ of a true Dirac operator on the underlying noncommutative space, see \cite{Skandalisgnc}. \\ 
We recall that a spectral triple $ (\A, \mathcal{H}, D) $ is called regular if for each $ a \in \A $ both $ a $ and $ [D, a] $ are contained in the 
domain of all powers of the derivation $ \delta $ given by 
$$ 
\delta(T) = [|D|, T], 
$$ 
see \cite{GFVbook}. 
\begin{lemma} \label{regularityofozawarieffelsptriple}
The above spectral triple $ (\A, \H, D) $ associated to a filtered algebra $ \A $ with a state is regular.
\end{lemma}
\proof Notice again that we have $ |D| = D $ in this case. In order to verify regularity we will closely follow the proof of 
lemma 1.1 of \cite{OzawaRieffel}. Let $ a $ be an element of $ \A_p $. Then we have 
$$ 
[D, a] = \sum_{j \mid |j| \leq p} j T_j 
$$ 
where $ T_j = \sum_m p_m a p_{m - j} $ is a bounded operator for each $ j $. 
We consider $ [D, [D, a]] = \sum_{j \mid |j| \leq p} [D, T_j] $. Now 
$$ 
[D, T_j] = [D, \sum_m p_m a p_{m - j}] = \sum_m j p_m a p_{m - j} = j T_j. 
$$
Thus, $ [D, [D, a]] = \sum_{j \mid |j| \leq p} j^2 T_j. $ The statement of the corollary follows by repeated application of this technique. \qed \\
Now we turn to the case of quantum groups. Let $ G $ be a discrete quantum group and $ l $ be a proper length function on it. Imitating 
Connes' construction in the group case, we construct a spectral triple with the Hilbert space $ l^2(G) $ where $ \mathbb{C}[G] $ acts by the regular representation. Our candidate for the Dirac operator is $ D $ where
$$ 
\dom(D) = \biggl\{\sum_{\alpha \in \Irr(G)} \sum_{i,j = 1}^{\dim(\alpha)} a^\alpha_{ij} \hat{\Lambda}(u^\alpha_{ij}) \mid 
\sum_{\alpha,i,j} l(\alpha)^2 |a^\alpha_{ij} |^2 \|\hat{\Lambda}(u^\alpha_{ij})\|^2 < \infty \biggr\}
$$ 
and 
$$ 
D\biggl(\sum_{\alpha, i, j} a^{\alpha}_{ij} \hat{\Lambda}(u^{\alpha}_{ij})\biggr) = \sum_{\alpha,i,j} l(\alpha) a^{\alpha}_{ij} \hat{\Lambda}(u^\alpha_{ij}). 
$$
The algebra $ \A = \mathbb{C}[G] $ is naturally a filtered algebra with the filtration defined by letting $ \A_n $ be the linear span of 
all matrix coefficients $ u^{\alpha}_{ij} $ such that $ l(\alpha) \leq n $. Since we assume that the length function $ l $ is proper there are only 
finitely many irreducible corepresentations of length $ n $. Hence the subspaces $ \A_n $ are all finite dimensional. Moreover
$ \A_0 = \mathbb{C}1 $, and it is easy 
to see from the properties of length functions that  $ \A^*_n = \A_n $ and $ \A_n \A_m \subset \A_{n + m} $.   
\begin{lemma} \label{regularspectraltriple}
Let $ G $ be a discrete quantum group with a proper length function. Then the triple $ (\mathbb{C}[G], l^2(G), D) $ constructed above is a 
regular spectral triple.
\end{lemma}
\proof Let $ A = C^*_\red(G) $ and $ \A = \mathbb{C}[G] $. The Hilbert space $ l^2(G) $ is the GNS-Hilbert space with respect to the Haar 
state $ \hat{\phi} $, and $ \hat{\phi} $ is  faithful on $ \mathbb{C}[G] $. By the observations made above, the 
filtration $ (\A_m)_{m \in \mathbb{N}_0} $ satisfies the conditions of lemma \ref{ozawarieffelsptriple}, and therefore $ (\mathbb{C}[G], l^2(G), D) $ 
is a spectral triple. Regularity follows from lemma \ref{regularityofozawarieffelsptriple}. \qed \\
Recall that a spectral triple $ (\A, \H, D) $ is called $ p $-summable if $ \Tr(\hat{D}^{-p}) < \infty $, where $ \hat{D} $ denotes the restriction 
of $ |D| $ to the orthogonal complement of its kernel. \\
The following result is a variant of proposition 6 in \cite{Connescms}. 
\begin{prop} \label{summability}
Let $ G $ be an amenable discrete quantum group of rapid decay, and fix a proper length function and constants $ c, s $ such that 
$$
\|a\|_{op} \leq c \|a\|_{2,s} 
$$ 
for all $ a \in \mathbb{C}[G] $. Then the associated spectral triple $ (\CH[G], l^2(G), D) $ constructed above is $ p $-summable for all $ p > 2s + 1 $. 
\end{prop}
\proof Since $ G $ is amenable, we see from the proof of proposition \ref{rdpgequiv} that 
$$ 
\sum_{l (\alpha) \in (n - 1, n]} \dim(\alpha)^2 = \phi(q_n) \leq c^2 (1 + n)^{2s}. 
$$ 
Using this, we find
\begin{align*}
\Tr(\hat{D}^{-p}) &= \sum_{\alpha \neq \epsilon} \dim(\alpha)^2 l(\alpha)^{-p} \\
&\leq r \sum_{l(\alpha) \in (0,1]} \dim(\alpha)^2 + \sum_{n > 1} \sum_{l(\alpha) \in (n - 1, n]} \dim(\alpha)^2 (n - 1)^{-p} \\
&\leq r \sum_{l(\alpha) \in (0,1]} \dim(\alpha)^2 + \sum_{n > 1} \frac{c^2 (n + 1)^{2s}}{(n - 1)^p} 
\end{align*}
where $ r $ is the inverse of the smallest nonzero value of $ l $ in the interval $ (0,1] $. This yields the claim. \qed \\
Notice that we cannot drop the amenability assumption in proposition \ref{summability} in general. For instance, according to 
\cite{Vergniouxrd}, the free orthogonal quantum group $ \mathbb{F}O(n) $ has property RD, but the corresponding spectral triple
will fail to be finitely summable as soon as $ n > 2 $.

\section{Compact quantum metric spaces} \label{seccompactquantummetricspace}

We recall that the geodesic distance on a compact Riemannian spin manifold can be recovered from the associated spectral triple. 
For a general spectral triple $ (\A, \H, D) $, one may take 
$$
d(\mu, \nu) = \sup_{a \in \A \mid \| [D, a] \| \leq 1}  |\mu(a) - \nu(a)| 
$$
as an Ansatz to define a metric on the state space $ S(A) $ of the $ C^* $-algebra closure $ A $ of $ \A $ in $ \LH(\H) $, 
thus generalising the Monge-Kantorovich metric on probability measures. 
We note that without further assumptions the above formula may yield $ d(\mu, \nu) = \infty $ for some states, 
see \cite{Connesgravmat}. \\ 
The above considerations, along with further examples, motivated Rieffel to introduce the concept of a quantum metric 
space. In \cite{Rieffelcqms} the theory is developed starting from order unit spaces instead of $ C^* $-algebras, but for 
our purposes the following definition is sufficient. 
\begin{definition} \label{defcqms}
Let $ A $ be a unital $ C^* $-algebra and let $ \A \subset A $ be a dense unital $ * $-subalgebra. A Lipschitz seminorm on $ \A $ is a seminorm 
$ L: \A \rightarrow [0, \infty) $ such that $ L(a^*) = L(a) $ for all $ a \in \A $ and $ L(a) = 0 $ iff $ a \in \mathbb{C} 1 $. 
A Lipschitz seminorm is called a Lip-norm if the topology on $ S(A) $ induced by 
$$
d_L(\mu, \nu) = \sup_{a \in \A \mid L(a) \leq 1}  |\mu(a) - \nu(a)| 
$$
coincides with the $ w^* $-topology. \\ 
A unital $ C^* $-algebra $ A $ together with a Lip-norm on a dense $ * $-subalgebra $ \A \subset A $ is called a compact quantum metric space.
\end{definition} 
Starting from definition \ref{defcqms} it is natural to ask for conditions ensuring that a Lipschitz seminorm is indeed a Lip-norm. 
Ozawa-Rieffel gave the following criterion, see proposition 1.3 in \cite{OzawaRieffel}. 
\begin{prop} \label{ozawarieffellipnorm}
Let $ A $ be a unital $ C^* $-algebra, and let $ L $ be a Lipschitz seminorm on a dense unital $ * $-subalgebra $ \A \subset A $. 
If $ \sigma $ is a state on $ A $ such that 
$$
E = \{a \in A \mid L(a) \leq 1 \; \text{and} \; \sigma(a) = 0 \} 
$$
is totally bounded with respect to the norm of $ A $, then $ L $ is a Lip-norm. 
\end{prop} 
Recall that a subset $ B $ of a metric space $ X $ is totally bounded if for any $ \epsilon > 0 $ there exists 
finitely many balls of radius $ \epsilon $ whose union covers $ B $. \\ 
Let us now assume that $ (\A, \H, D) $ is the spectral triple associated to a filtered algebra as in section \ref{secspectraltriples}. One 
may ask whether the resulting Lipschitz seminorm $ L(a) = \|\delta(a)\| = \|[D,a]\| $ is a Lip-norm on $ \A $, viewed as a dense subalgebra 
of its norm closure $ A $. This seems to be unclear in general. \\ 
In \cite{OzawaRieffel}, Ozawa and Rieffel showed that $ L $ is indeed a Lip-norm provided a condition of the form
$$ 
\| p_m a p_n   \| \leq c \| a \|_2 
$$ 
holds for $ a \in \A_k $ and all $ k,m,n $. They call such an inequality a Haagerup type condition. As remarked in section \ref{secspectraltriples}, 
Connes' spectral triples on group algebras coming from length functions are special cases of the Ozawa-Rieffel construction in 
lemma \ref{ozawarieffelsptriple}. Examples of groups for which the Haagerup type condition hold include word hyperbolic 
groups \cite{Connesbook}, \cite{OzawaRieffel} as well as free products of the form $ G_1 \ast G_2 $ where $ G_1 $ and $ G_2 $ satisfy the 
Haagerup type condition and one works with tracial states on the group algebras of $ G_1 $ and $ G_2 $. \\ 
It is not clear, however, whether $ L $ is a Lip-norm for groups of rapid decay in general. In this case one has a weaker inequality of the form
$$ 
\| p_m a p_n \| \leq c P(k) \| a \|_2 
$$
for all $ a \in \A_k $, where $ P $ is a polynomial. Antonescu and Christensen observed that one can easily prove 
the Lip-norm property in this case if one works with a slightly different Lipschitz seminorm instead \cite{ACmetricsgroup}. More precisely, 
for $ k \in \mathbb{N} $ and $ a \in \A = \mathbb{C}[G] $, consider
$$ 
L^k(a) = \| \delta^k(a) \| = \|[D, [D, \dots [D,a] \dots]\|, 
$$ 
where the commutator is taken $ k $ times. \\
We will consider a similar construction in the case of discrete quantum groups, and prove some lemmas which are needed later on. However, in order 
to accommodate non-unimodular discrete quantum groups, we need to use a twisted version of the seminorm $ L^k $.    
Throughout, we assume that $ G $ is a discrete quantum group equipped with a proper length function $ l $. 
\begin{lemma} \label{lemma1}
Using the same notation as before, we have 
$ \delta^k(u^{\alpha}_{ij}) \hat{\Lambda}(1) = D^k \hat{\Lambda}(u^{\alpha}_{ij}) $ for all $ k \in \mathbb{N} $. 
In particular, for a finite sum of the form 
$ a = \sum_{\alpha,i,j} a^{\alpha}_{ij} u^{\alpha}_{ij} $, we have 
$ \delta^k(a) \hat{\Lambda}(1) = \sum_{\alpha,i,j} l(\alpha)^k a^{\alpha}_{ij} \hat{\Lambda}(u^{\alpha}_{ij}) $.  
\end{lemma}
\proof We proceed by induction. For $ k = 1 $, we have
\begin{align*} 
\delta(u^{\alpha}_{ij}) \hat{\Lambda}(1) &= [D, u^{\alpha}_{ij}] \hat{\Lambda}(1) \\
&= D \hat{\Lambda}(u^{\alpha}_{ij}) - u^{\alpha}_{ij} D\hat{\Lambda}(1) = D\hat{\Lambda}(u^{\alpha}_{ij}). 
\end{align*}
Let us now assume that $ \delta^k(u^{\alpha}_{ij}) \hat{\Lambda}(1) = D^k \hat{\Lambda}(u^{\alpha}_{ij}) $. Then 
\begin{align*} 
\delta^{k + 1}(u^{\alpha}_{ij}) \hat{\Lambda}(1) &= [D, \delta^k(u^{\alpha}_{ij})] \hat{\Lambda}(1) \\
&= D \delta^k (u^{\alpha}_{ij}) \hat{\Lambda}(1) - \delta^k(u^{\alpha}_{ij}) D \hat{\Lambda}(1) \\
&= D D^k \hat{\Lambda}(u^{\alpha}_{ij}) \\
&= D^{k + 1} \hat{\Lambda}(u^{\alpha}_{ij}), 
\end{align*}
which yields the claim. \qed 
\begin{lemma} \label{lemma2}
If $ a \in \mathbb{C}[G] $ is such that $ L^k(a) = 0 $, then $ a $ is a scalar multiple of the identity. Hence $ L^k $ is a 
Lipschitz seminorm.
\end{lemma}
\proof Let $ a = \sum_{\alpha, i,j} a^{\alpha}_{ij} u^{\alpha}_{ij} \in \CH[G] $. By definition, the relation $ L^k(a) = 0 $ 
implies that $ \delta^k(a) = 0 $. In particular we have $ \delta^k(a) \hat{\Lambda}(1) = 0 $. Therefore, according to lemma \ref{lemma1} we obtain
$$ 
0 = \sum_{\alpha,i,j} a^{\alpha}_{ij} D^k \hat{\Lambda}(u^{\alpha}_{ij}) 
= \sum_{\alpha,i,j} l(\alpha)^k a^{\alpha}_{ij}  \hat{\Lambda}(u^{\alpha}_{ij}) = 0. 
$$
Since the vectors $ \hat{\Lambda}(u^{\alpha}_{ij}) $ form a linearly independent set in $ l^2(G) $, we conclude $ l(\alpha)^k a^{\alpha}_{ij} = 0 $ 
for all $ \alpha,i,j $. If $ l(\alpha) \neq 0 $ then $ a^{\alpha}_{ij} = 0 $. Since $ l $ is a proper length function we conclude that $ a $ has 
to be a scalar multiple of the identity. Moreover, we clearly have $ L^k(1) = 0 $ and $ L^k(a) = L^k(a^*) $. 
Hence $ L^k $ is a Lipschitz seminorm. \qed \\
Fix a natural number $ k $ and consider the operator $ T = C^{\frac{1}{2k}} $. As before, we let $ \delta(a) = [D,a] $ and define 
$$ 
\delta_{T} (a) = [D, T a T], \qquad L^k_{T} (a) = \| \delta^k_T(a) \|
$$
for $ a \in \A = \CH[G] $. 
\begin{lemma}
Let $ k \in \mathbb{N} $. For any $ a \in \A $ we have $ \delta^k_T(a) = T^k \delta^k(a) T^k $, and $ L^k_T $ is a Lipschitz seminorm on $ \A $.
\end{lemma}
\proof Since $ T $ commutes with $ D $ we have $ \delta_T(a) = [D, T a T] = T [D,a] T $, and so $ \delta^k_T(a) = T^k \delta^k(a) T^k $. 
Moreover 
$$ 
L^k_T(a^*) = \|\delta_T^k(a^*)\| = \| T^k \delta^k(a^*) T^k \| = \| T^k \delta^k(a) T^k \| = L^k_T(a) 
$$
since $ T $ is self-adjoint. For the Lipschitz seminorm property, we need to check that if $ L^k_T(a) = 0 $, then $ a $ is a scalar multiple of 
the identity. But $ L^k_T(a) = 0 $ implies that $ T^k \delta^k(a) T^k = 0 $. Since $ T $ is invertible we conclude $ \delta^k(a) = 0 $ 
which means $ L^k(a) = 0 $. Hence the desired conclusion follows from lemma \ref{lemma2}. \qed 
\begin{lemma} \label{labelledon21stapril}
Let $ k \in \mathbb{N} $. Then we have 
\begin{align*}
\sum_{\alpha, i, j} \dim(\alpha)^{-1} l(\alpha)^{2k} | a^{\alpha}_{ij}|^2
&= \|\sum_{\alpha, i, j} \delta^{k}_T (a^{\alpha}_{ij} u^\alpha_{ij}) \hat{\Lambda}(1) \|^2_2 \\
&\leq L^k_T\biggl(\sum_{\alpha,i,j} a^{\alpha}_{ij} u^\alpha_{ij}\biggr)^2.
\end{align*}
\end{lemma}
\proof Using lemma \ref{lemma1} we compute
\begin{align*}
&\|\sum_{\alpha, i, j} \delta^{k}_T (a^{\alpha}_{ij} u^\alpha_{ij}) \hat{\Lambda}(1) \|^2_2 
= \| \sum_{\alpha,i,j} a^{\alpha}_{ij} l(\alpha)^k T^k \hat{\Lambda}(u^{\alpha}_{ij}) \|^2_2 \\
&= \sum_{\alpha} \sum_{i,j} \sum_{i^{\prime},j^{\prime}} 
\overline{a^{\alpha}_{ij}} a^{\alpha}_{i^{\prime}j^{\prime}} \frac{\dim_q(\alpha)}{\dim(\alpha)} 
l(\alpha)^{2k} (F^{\alpha})^{\frac{1}{2}}_{ii} (F^{\alpha})^{\frac{1}{2}}_{i^{\prime}i^{\prime}} 
\bra \hat{\Lambda}(u^{\alpha}_{ij}), \hat{\Lambda}(u^{\alpha}_{i^{\prime} j^{\prime}}) \ket \\
&= \sum_{\alpha,i,j} \dim(\alpha)^{-1} l(\alpha)^{2k} | a^\alpha_{ij}|^2. 
\end{align*}
This yields the claim. \qed

\section{Rapid decay and the Lip norm property} \label{seclipnorm}

In this section we show that the Lipschitz norms considered above have the Lip norm property for all sufficiently large exponents 
provided the quantum group under consideration has property RD. \\ 
For simplicity we concentrate on the case of a finitely generated discrete quantum group $ G $ equipped with a word length function $ l $. 
Moreover we assume that $ G $ has property RD, and we fix constants $ c, s > 0 $ such that 
$$ 
\| a \|_{op} \leq c \| a \|_{2,s} 
$$ 
holds for all $ a $ in $ \CH[G] $ with respect to the Sobolev norms associated with $ l $. \\
Using the notation introduced in section \ref{seccompactquantummetricspace}, we will apply proposition \ref{ozawarieffellipnorm} by Ozawa-Rieffel to show 
that there exists a positive integer $ k $ such that $ L^k_T $ has the Lip norm property. We start with the following lemma.
\begin{lemma} \label{lipnorm}
Let $ a = \sum_{\alpha,i,j} a^{\alpha}_{ij} u^\alpha_{ij} \in \CH[G] $. Then we have
$$
\|a\|^2_{2,s} = \sum_{\alpha, i,j} \frac{1}{\dim(\alpha)} (1 + l(\alpha))^{2s} |a^\alpha_{ij}|^2. 
$$ 
\end{lemma}
\proof We compute 
\begin{align*} 
\| a \|^2_{2,s} &= \sum_{\alpha, i,j} 
\bra (1 + L)^s a^\alpha_{ij} \F^{-1}(u^\alpha_{ij}) C^{1/2}, (1 + L)^s a^\alpha_{ij} \F^{-1}(u^\alpha_{ij}) C^{1/2} \ket \\
&= \sum_{\alpha, i,j} \frac{\dim_q(\alpha)}{\dim(\alpha)} (1 + l(\alpha))^{2s} |a^\alpha_{ij}|^2 (F^\alpha)_{ii}
\bra u^\alpha_{ij}, u^\alpha_{ij} \ket \\
&= \sum_{\alpha, i,j} \frac{1}{\dim(\alpha)} (1 + l(\alpha))^{2s} |a^\alpha_{ij}|^2. 
\end{align*}
This yields the claim. \qed \\
Now, for a positive integer $ k > s $, let 
$$ 
E = \{a \in \mathbb{C}[G] \mid L^{k}_T(a) \leq 1 \; \text{and} \; \hat{\phi}(a) = 0 \}, 
$$
which we view as a subset of $ \CH[G] \subset C^*_\red(G) $. 
\begin{lemma} \label{termE1}
For all $ n \in \mathbb{N} $ there exists a constant $ c_n $ such that 
$$ 
\| \sum_{\alpha, i, j \mid l(\alpha) \leq n} a^{\alpha}_{ij} u^{\alpha}_{ij}  \|_{op} \leq c_n 
$$ 
for any $ a = \sum_{\alpha, i, j} a^{\alpha}_{ij} u^{\alpha}_{ij} \in E $. 
\end{lemma}
\proof We recall that $ \epsilon $ denotes the trivial corepresentation. Using $ a^\epsilon_{11} = \hat{\phi}(a) = 0 $ and the Cauchy-Schwarz inequality we compute
\begin{align*} 
&\sum_{\alpha, i, j \mid l(\alpha) \leq n} |a^{\alpha}_{ij}|
= \sum_{\alpha, i, j \mid 1 \leq l(\alpha) \leq n} l(\alpha)^{-k} 
\sqrt{\dim(\alpha)} l(\alpha)^k |a^{\alpha}_{ij}| 
\sqrt{\dim(\alpha)}^{\,-1} \\
&\leq \biggl(\sum_{\alpha, i, j \mid l(\alpha) \leq n}  l(\alpha)^{-2k}  \dim(\alpha) \biggr)^{\frac{1}{2}} 
\biggl(\sum_{\alpha, i, j \mid l(\alpha) \leq n} l(\alpha)^{2k} |a^{\alpha}_{ij}|^2 \dim(\alpha)^{-1} \biggr)^{\frac{1}{2}} \\
&\leq c_n L_T^k(a) \leq c_n,  
\end{align*}
where we choose $ c_n $ such that 
$$ 
\biggl(\sum_{\alpha, i, j \mid l(\alpha) \leq n}  l(\alpha)^{-2k}  \dim(\alpha) \biggr)^{\frac{1}{2}} \leq c_n, 
$$ 
and we use lemma \ref{labelledon21stapril}. \\
Since the matrices $ (u^\alpha_{ij}) $ are unitary we have $ \|u^{\alpha}_{ij}\|_{op} \leq 1 $ for all $ \alpha $ and $ i,j = 1, \dots, \dim(\alpha) $. 
Hence we obtain 
\begin{align*} 
\| \sum_{\alpha, i, j \mid l(\alpha) \leq n} a^{\alpha}_{ij} u^{\alpha}_{ij} \|_{op} 
&\leq \sum_{\alpha, i, j \mid l(\alpha) \leq n} |a^{\alpha}_{ij}| \|u^{\alpha}_{ij}\|_{op}  
\leq \sum_{\alpha, i, j \mid l(\alpha) \leq n} |a^{\alpha}_{ij}| \leq c_n.  
\end{align*} 
This yields the claim. \qed 
\begin{lemma} \label{termE2}
Let $ a = \sum_{\alpha,i,j} a^{\alpha}_{ij} u^\alpha_{ij} \in \CH[G] $ and fix a positive integer $ k > s $. Then we have
$$
\|\sum_{\alpha,i,j, l(\alpha) > n} a^{\alpha}_{ij} u^\alpha_{ij} \|^2_{op} 
\leq c^2 2^{2s} n^{2(s - k)} \sum_{\alpha, i, j, l(\alpha) > n} \dim(\alpha)^{-1} l(\alpha)^{2 k} |a^{\alpha}_{ij}|^2
$$
for any $ n \in \mathbb{N} $. 
\end{lemma}
\proof We follow Antonescu-Christensen \cite{ACmetricsgroup}. For $ \alpha $ satisfying $ l(\alpha) > n $ write 
$$ 
(1 + l(\alpha))^{2s} \leq 2^{2s} l(\alpha)^{2s} \leq 2^{2s} n^{2s - 2 k} l(\alpha)^{2 k}.
$$ 
Therefore, by property RD and lemma \ref{lipnorm}, we have
\begin{align*} 
&\|\sum_{\alpha,i,j, l(\alpha) > n} a^{\alpha}_{ij} u^\alpha_{ij} \|^2_{op}
\leq c^2 \| \sum_{\alpha,i,j, l(\alpha) > n} a^{\alpha}_{ij} u^\alpha_{ij} \|^2_{2,s} \\
&= c^2 \sum_{\alpha,i,j, l(\alpha) > n} \dim(\alpha)^{-1} (1 + l(\alpha))^{2s} |a^{\alpha}_{ij}|^2 \\
&\leq c^2 2^{2s} n^{2(s - k)} \sum_{\alpha, i, j, l(\alpha) > n} \dim(\alpha)^{-1} l(\alpha)^{2 k} |a^{\alpha}_{ij}|^2.
\end{align*}
This finishes the proof. \qed \\
We are now ready to prove the following result. 
\begin{theorem}
For any positive integer $ k > s $, the algebra $ \CH[G] \subset C^*_\red(G) $ together with the seminorm $ L^{k}_T $ is a compact quantum metric 
space.
\end{theorem}
\proof We have to show that $ L^{k}_T $ is a Lip norm. Fix $ \epsilon > 0 $ and choose $ n $ such that $ c 2^{s} n^{s - k} < \epsilon $ 
and write $ E = E_1 + E_2, $ where 
$$ 
E_1 = \{x \in E \mid x = \sum_{\alpha, i, j \mid l (\alpha) \leq n} a^{\alpha}_{ij} u^{\alpha}_{ij} \} 
$$ 
and 
$$ 
E_2 = \{x \in E \mid x = \sum_{\alpha,i,j \mid l(\alpha) > n} a^{\alpha}_{ij} u^{\alpha}_{ij} \}. 
$$
Using lemma \ref{termE1} we see that $ E_1 $ is a bounded subset of a finite dimensional normed space, and thus totally bounded. 
Moreover, by our choice of $ n $ and lemma \ref{termE2}, the set $ E_2 $ is contained in the $ \epsilon $-ball around $ 0 $ in $ C^*_\red(G) $. 
This completes the proof. \qed

\bibliographystyle{plain}

\bibliography{cvoigt}

\def\cprime{$'$} \def\cprime{$'$} \def\cprime{$'$} \def\cprime{$'$}
  \def\polhk#1{\setbox0=\hbox{#1}{\ooalign{\hidewidth
  \lower1.5ex\hbox{`}\hidewidth\crcr\unhbox0}}} \def\cprime{$'$}
  \def\cprime{$'$} \def\Dbar{\leavevmode\lower.6ex\hbox to 0pt{\hskip-.23ex
  \accent"16\hss}D} \def\cftil#1{\ifmmode\setbox7\hbox{$\accent"5E#1$}\else
  \setbox7\hbox{\accent"5E#1}\penalty 10000\relax\fi\raise 1\ht7
  \hbox{\lower1.15ex\hbox to 1\wd7{\hss\accent"7E\hss}}\penalty 10000
  \hskip-1\wd7\penalty 10000\box7}
  \def\cfudot#1{\ifmmode\setbox7\hbox{$\accent"5E#1$}\else
  \setbox7\hbox{\accent"5E#1}\penalty 10000\relax\fi\raise 1\ht7
  \hbox{\raise.1ex\hbox to 1\wd7{\hss.\hss}}\penalty 10000 \hskip-1\wd7\penalty
  10000\box7}
\begin{thebibliography}{10}

\bibitem{ACmetricsgroup}
Cristina Antonescu and Erik Christensen.
\newblock Metrics on group {$C^*$}-algebras and a non-commutative
  {A}rzel\`a-{A}scoli theorem.
\newblock {\em J. Funct. Anal.}, 214(2):247--259, 2004.

\bibitem{BSUM}
Saad Baaj and Georges Skandalis.
\newblock Unitaires multiplicatifs et dualit\'e pour les produits crois\'es de
  {$C^*$}-alg\`ebres.
\newblock {\em Ann. Sci. \'Ecole Norm. Sup. (4)}, 26(4):425--488, 1993.

\bibitem{Banicacompactsubfactor}
Teodor Banica.
\newblock Representations of compact quantum groups and subfactors.
\newblock {\em J. Reine Angew. Math.}, 509:167--198, 1999.

\bibitem{CPspectraltriple}
Partha~Sarathi Chakraborty and Arupkumar Pal.
\newblock Equivariant spectral triples on the quantum {${\rm SU}(2)$} group.
\newblock {\em $K$-Theory}, 28(2):107--126, 2003.

\bibitem{Connescms}
A.~Connes.
\newblock Compact metric spaces, {F}redholm modules, and hyperfiniteness.
\newblock {\em Ergodic Theory Dynam. Systems}, 9(2):207--220, 1989.

\bibitem{Connesbook}
Alain Connes.
\newblock {\em Noncommutative geometry}.
\newblock Academic Press Inc., San Diego, CA, 1994.

\bibitem{Connesgravmat}
Alain Connes.
\newblock Gravity coupled with matter and the foundation of non-commutative
  geometry.
\newblock {\em Comm. Math. Phys.}, 182(1):155--176, 1996.

\bibitem{DLSSVdirac}
Ludwik D{\c{a}}browski, Giovanni Landi, Andrzej Sitarz, Walter van Suijlekom,
  and Joseph~C. V{\'a}rilly.
\newblock The {D}irac operator on {${\rm SU}_q(2)$}.
\newblock {\em Comm. Math. Phys.}, 259(3):729--759, 2005.

\bibitem{GFVbook}
Jos{\'e}~M. Gracia-Bond{\'{\i}}a, Joseph~C. V{\'a}rilly, and H{\'e}ctor
  Figueroa.
\newblock {\em Elements of noncommutative geometry}.
\newblock Birkh\"auser Advanced Texts: Basler Lehrb\"ucher. [Birkh\"auser
  Advanced Texts: Basel Textbooks]. Birkh\"auser Boston Inc., Boston, MA, 2001.

\bibitem{Haagerupnonnuclear}
Uffe Haagerup.
\newblock An example of a nonnuclear {$C^{\ast} $}-algebra, which has the
  metric approximation property.
\newblock {\em Invent. Math.}, 50(3):279--293, 1978/79.

\bibitem{Jolissaintrd}
Paul Jolissaint.
\newblock Rapidly decreasing functions in reduced {$C^*$}-algebras of groups.
\newblock {\em Trans. Amer. Math. Soc.}, 317(1):167--196, 1990.

\bibitem{KVLCQG}
Johan Kustermans and Stefaan Vaes.
\newblock Locally compact quantum groups.
\newblock {\em Ann. Sci. \'Ecole Norm. Sup. (4)}, 33(6):837--934, 2000.

\bibitem{NTDirac}
Sergey Neshveyev and Lars Tuset.
\newblock The {D}irac operator on compact quantum groups.
\newblock {\em J. Reine Angew. Math.}, 641:1--20, 2010.

\bibitem{NVpoincare}
Ryszard Nest and Christian Voigt.
\newblock Equivariant {P}oincar\'e duality for quantum group actions.
\newblock {\em J. Funct. Anal.}, 258(5):1466--1503, 2010.

\bibitem{OzawaRieffel}
Narutaka Ozawa and Marc~A. Rieffel.
\newblock Hyperbolic group {$C^*$}-algebras and free-product {$C^*$}-algebras
  as compact quantum metric spaces.
\newblock {\em Canad. J. Math.}, 57(5):1056--1079, 2005.

\bibitem{Rieffelcqms}
Marc~A. Rieffel.
\newblock Compact quantum metric spaces.
\newblock In {\em Operator algebras, quantization, and noncommutative
  geometry}, volume 365 of {\em Contemp. Math.}, pages 315--330. Amer. Math.
  Soc., Providence, RI, 2004.

\bibitem{Skandalisgnc}
G.~Skandalis.
\newblock G\'eom\'etrie non commutative d'apr\`es {A}lain {C}onnes: la notion
  de triplet spectral.
\newblock {\em Gaz. Math.}, (94):44--51, 2002.

\bibitem{vDadvances}
A.~Van~Daele.
\newblock An algebraic framework for group duality.
\newblock {\em Adv. Math.}, 140(2):323--366, 1998.

\bibitem{Vergniouxrd}
Roland Vergnioux.
\newblock The property of rapid decay for discrete quantum groups.
\newblock {\em J. Operator Theory}, 57(2):303--324, 2007.

\bibitem{Voiculescuquasicentral}
Dan Voiculescu.
\newblock On the existence of quasicentral approximate units relative to normed
  ideals. {I}.
\newblock {\em J. Funct. Anal.}, 91(1):1--36, 1990.

\bibitem{Woronowiczleshouches}
S.~L. Woronowicz.
\newblock Compact quantum groups.
\newblock In {\em Sym\'etries quantiques ({L}es {H}ouches, 1995)}, pages
  845--884. North-Holland, Amsterdam, 1998.

\end{thebibliography}

\end{document}